\documentclass[12pt]{amsart}

\usepackage{amssymb}
\usepackage{psfrag}

\newtheorem{thm}{Theorem}%[section]

\newtheorem{lem}[thm]{Lemma}
\newtheorem{cor}[thm]{Corollary}

\newtheorem{prop}[thm]{Proposition}

\newtheorem{conj}[thm]{Conjecture}
\theoremstyle{definition}
\newtheorem{defn}[thm]{Definition}

\newtheorem{say}[thm]{}

\newtheorem{ques}[thm]{Question}    %!!!!!!!!!!!!!!!!!!!!

\newtheorem{rem}[thm]{Remark}          
\newtheorem{ack}{Acknowledgments}

\newtheorem{defn-thm}[thm]{Definition--Theorem}  %!!!!!!!!!!!!!!!!!!!!!!!!
\theoremstyle{remark}
\newtheorem{claim}[thm]{Claim}

\renewcommand{\c}[0]{{\mathbb C}}  

\renewcommand{\o}[0]{{\mathcal O}} 

\newcommand{\n}[0]{{\mathbb N}}
\renewcommand{\r}[0]{{\mathbb R}} 

\renewcommand{\a}[0]{{\mathbb A}}

\newcommand{\p}[0]{{\mathbb P}}
\newcommand{\f}[0]{{\mathbb F}}
\newcommand{\q}[0]{{\mathbb Q}}

\newcommand{\qtq}[1]{\quad\mbox{#1}\quad}
\newcommand{\spec}[0]{\operatorname{Spec}}

\newcommand{\gal}[0]{\operatorname{Gal}}

\newcommand{\Hom}[0]{\operatorname{Hom}}

\newcommand{\hilb}[0]{\operatorname{Hilb}}

\newcommand{\CH}[0]{\operatorname{CH}}

\def\into{\DOTSB\lhook\joinrel\rightarrow}
\def\to {\longrightarrow}

\begin{document}
\bibliographystyle{amsplain}

\title[Rationally connected varieties]
{Rationally connected varieties
over finite  fields}
\author{J\'anos Koll\'ar and Endre Szab\'o}

\maketitle

In this paper we study 
  rationally connected varieties 
defined over  finite fields. Then we lift these results to
rationally connected varieties 
 over local fields.

Roughly speaking, a variety $X$ over an algebraically closed field
is rationally connected if it contains a rational curve
through any number of assigned points $P_1,\dots,P_n$.
See \cite{Ko01} for an introduction to their theory
and for an explanation of the place of rationally
connected varieties in the framework of 
the classification theory 
of algebraic varieties.
There are many equivalent conditions defining this notion,
see, for instance, \cite[IV.3]{rcbook}.
The definition given  below essentially  corresponds to the case of 2 points.
In positive characteristic we also have to be mindful
of some inseparability problems.

\begin{defn}\label{SRC-defn}  A smooth, 
proper variety $X$ is called {\it separably rationally connected}
  or {\it SRC}, if there is a variety $U$ and a
morphism $F:U\times\p^1\to X$ such that the induced map
$$
F(\_\ ,(0{:}1))\times F(\_\ ,(1{:}0)) : U\to X\times X
%F\vert_{U\times\{0\}}\times F|_{U\times\{1\}} : U\to X\times X
$$
is dominant and separable.
\end{defn}

If $X$ is defined over a field $k$
and the $P_i$ are also defined over $k$, it is of interest
to find rational curves defined over $k$ passing through these points.
In general this cannot be done.
For instance, consider the surface
$$
S^0:=(x^2+y^2=-\prod_{i=1}^{2m} (z-a_i))\subset\a^3,\qtq{where $a_i\in \r$.}
$$
 Let $S\supset S^0$ be any smooth compactification. 
The fibers of the projection to the $z$-axis $S^0\to \a_z$
are conics, thus
$S$ is ruled over $\p^1$. Therefore  $S$ is rational over $\c$.
As a topological space $S(\r)$ has $m$ connected components. 
If we pick points
$P_1,P_2$ in different connected components of $S(\r)$
then there is no rational curve defined over $\r$
passing through both points (since $\p^1(\r)$ is connected).

The main technical  result of the present paper (Theorem  
\ref{main.geom.thm}) implies  that
for finite fields one can  almost always find such rational curves.
More precisely, we have the following:

\begin{thm}\label{main.charp.thm} There is a function $\Phi:\n^3\to\n$
with the following property: 

Let $X\subset \p^N$ be a smooth, projective,
separably rationally connected variety 
of dimension $\geq 3$
over a finite field $k$
and $S\subset X$ a zero--dimensional smooth subscheme. 

If $|k|>\Phi(\deg X ,\dim X ,\deg S)$, 
then there is a smooth  rational curve
$C_S\subset X$ containing $S$, where $|k|$ denotes the
cardinality of $k$.
\end {thm}

\begin{rem} (1) If $\dim X=2$ then there are very few smooth rational curves
even over $\bar k$. For instance, in $\p^2$ the only
smooth rational curves are lines and conics. If $S$
consists of 3 points on a line plus another point
then there is no smooth rational curve containing $S$.

For all applications considered in this paper,
one can work with $X\times \p^1$ instead of the surface $X$
to obtain the relevant conclusions.

One can also see from the proof that if $X$ is a surface,
we obtain an immersed curve $\p^1\to X$ through $S$.

(2)  We can also be more precise about the curve $C_S$.
For instance, with small modifications of the proof  one
sees that
 Theorem \ref{main.charp.thm} remains true if
we only assume that $S$ has embedding dimension 1.

We can also require that $C_S$ avoid a given codimension 2
subset of $X$ and that it be transversal
to a given divisor
at all points outside $S$.

(3) The function $\Phi$ is explicitly computable from the proof
though it is rather huge.
Even if we allow singular curves $C_S$, some condition
on $|k|$ is needed. For instance, the number of $k$-points in $S$
cannot exceed $|k|+1$, the number of $k$-points in $\p^1$.
We do not know any other necessary condition  with singular $C_S$.
In the smooth case, there are also problems 
 with too few $k$-points.
For instance, \cite[5.5]{SD} contains an example of a cubic
surface $S$ over $\f_2$  with only one rational point.
Thus $S\times S$ does  not contain any smooth rational 
curves defined over $\f_2$.
We  do not know what a reasonable $\Phi$  should be.

(4) By a recent result of \cite{esn}, an SRC variety defined
over a finite field always has a point in that field.

\end{rem}

We give two  applications of this result  to the
 study of the R-equivalence of points
and the Chow group of 0-cycles $\CH_0(X)$. The interesting part is
the 
Chow group of degree zero 0-cycles, denoted by  $\CH_0^0(X)$.

The Chow group of zero cycles
over a finite field is completely described in \cite{ka-sa}.
Their result implies
 that if $X$ is a smooth, projective
 variety 
defined over a finite field which is geometrically simply connected
then  $\CH_0^0(X)=0$. 
SRC varieties were recently shown to be  geometrically simply connected
(see \cite{kol-let} or \cite[3.6]{deb-b}),
so $\CH_0^0(X)=0$ for them.
 We give a geometric proof of this 
result for SRC varieties. 
The advantage of
this proof is that it allows one to go from finite fields
to local fields.
(For our purposes {\it local fields} are the quotient fields of
complete discrete valuation rings with finite residue field.) 
Lifting results
of this type have been known for
cubic surfaces \cite{SD}.
Recently this was
extended to  cubic hypersurfaces 
 \cite{ct-stud??}. 
For surfaces both results were established by \cite{ct-surfcase}.

\begin{thm}\label{finite-CH.thm}\cite{ka-sa} Let $k$ be a finite field
and $X$ a smooth, projective, separably rationally connected variety
over $k$. 
Then $\CH_0^0(X)=0$.
\end {thm}

\begin{thm}\label{local-CH.thm} Let $K$ be a local field 
with residue field $k$
and $X$ a variety over $K$. Assume that $X$ has a
smooth, projective, separably rationally connected reduction over $k$. 
Then $\CH_0^0(X)=0$.
\end {thm}

%%%%%%%%%%%%%%%%%%%%%%%%%%%%%%%%%%%%%%%%%

R-equivalence on cubic surfaces was introduced in \cite{manin}
in order to study the parametrization of rational points.
\cite{SD} proves that R-equivalence is trivial 
on  smooth cubic surfaces over finite fields. 
 Precise computations  on certain algebraic groups
and related  spaces are done in 
\cite{ gille1, gille2, mb01, mb02}.

\begin{defn}\label{R-defn}
Let $X$ be a proper 
variety defined over a field $k$. Two points  $x,x'\in X(k)$
are called directly R-equivalent if there is a morphism
$f:\p^1\to X$ (defined over $k$) such that $f(0{:}1)=x$
and $f(1{:}0)=x'$. The  equivalence relation  generated by this
is called {\it R-equivalence}. 
Informally speaking, two points are 
 {\it R-equivalent} if they can be connected by a chain
of rational curves over $k$. If all the points in $X(k)$ are
R-equivalent, then we say that  {\it R-equivalence is trivial}
on $X(k)$.
\end{defn}

\begin{thm}\label{finite-R.thm} There is a function $\Psi:\n^2\to\n$
with the following property: 

Let $X\subset \p^N$ be a smooth, projective,
separably rationally connected variety over the finite field $k$. 
If $|k|>\Psi(\deg X ,\dim X )$ then  R-equivalence  is trivial on $X(k)$.
\end {thm}

\begin{thm}\label{local-R.thm} There is a function $\Psi':\n^2\to\n$
with the following property: 

Let $K$ be a local field with
residue-field $k$ and $X\subset\p^N$ a variety over $K$. Assume that $X$ has a
smooth, projective, separably rationally connected reduction over $k$. 
If $|k|> \Psi'(\deg X ,\dim X )$ then  R-equivalence is trivial on $X(K)$.
\end {thm}

The following  corollary answers a question of Colliot-Th\'el\`ene
which 
was the starting point of our investigations.
 Let $X$ be a smooth projective  variety over a number
field $K$. For  a prime $P$, let $X_P$ denote the
corresponding variety over the $P$-adic completion 
$K_P$. Our results imply
that the Chow group of zero cycles and R-equivalence
are both trivial for almost every completion.
This has been known in several special cases,
see \cite{cts,  ct-surfcase, ctssd, gille2, mb01, mb02}.

\begin{cor} Let $X$ be a smooth projective 
rationally connected variety over a number
field $K$.  Then 
$\CH_0^0(X_P)=0$ and R-equivalence is trivial on $X(K_P)$
for  all but finitely many  primes $P$. 
\qed
\end{cor}

\begin{say}[Open problems]{\ }
There are several unsolved problems closely related to
these results. First of all, one may wish to know 
something about the function $\Psi(\deg X ,\dim X )$.
Unlike in the case of $\Phi$, the known examples suggest that
here the bound is not needed at all:

\begin{ques}\label{no-bounds.rem}
Do Theorems (\ref{finite-R.thm}) and  (\ref{local-R.thm})    remain true
without the assumption $|k|>\Psi(\deg X ,\dim X )$?
\end{ques}

More optimistically, one can ask the following 
 question, which could   be the
strongest version of (\ref{main.charp.thm}) that does not
need any assumption on $|k|$.

\begin{ques}\label{main.thm.ques}
 Let $X$ be a smooth,
projective, separably rationally connected variety over a finite field $k$. 
Let  $C$ be a  smooth curve over $k$ and $T\subset C$ 
 a zero--dimensional subscheme. 
Is it true that every morphism
$g_T:T\to X$   can be extended to $g:C\to X$?
\end{ques}

This is easily seen to be true if $X$ is rational over $k$.

Our current methods
say very little about
SRC varieties over local fields with bad reduction.
The following basic finiteness result for R-equivalence is
known:

\begin{thm}\label{local-R-veges.thm}\cite{rcloc}
 Let $X$ be a smooth,
projective, separably rationally connected variety over a local field
$K$. 
Then $X(K)$ has finitely many R-equivalence classes.
\end {thm}

The analogous statement about the Chow group is not known,
except for surfaces \cite{ct-surfcase}.

\begin{conj}\label{local-CH-veges.thm} Let $X$ be a smooth,
projective, separably rationally connected variety defined over a
local field.  Then  $\CH_0^0(X)$  is finite.
\end {conj}

The separability condition in (\ref{local-CH.thm}) is essential.
Indeed, there are varieties of general type 
over $\q_p$ whose reduction to $\f_p$ is smooth and
purely inseparably uniruled.
(Some very nice examples are given in \cite{shi-ka}.)
The Chow group  $\CH_0^0$ is
very large for some of them, see \cite{mumf}.

On the other hand, one would expect that our results 
hold for Fano
varieties even if they are not separably rationally connected.
The first case to investigate would be
smooth hypersurfaces $X\subset \p^n$ of degree $\leq n$.
These are all rationally connected but not 
known to be always SRC in low characteristics.
(The examples in \cite[V.5]{rcbook} have some mild singularities.) 
\end{say}

\begin{say}
Let $X$ be a smooth, projective,
separably rationally connected variety 
over a   field $k$ and 
 $S\subset X$  a zero--dimensional  smooth  subscheme.
Our aim is  to find a rational curve $C\subset  X$ 
defined over $k$, passing through
$S$. 

If $X$ is SRC, there are  such curves over $\bar k$
(cf. \cite[IV.3.9]{rcbook}),
thus the space $M\subset \hilb(X)$
of  smooth rational curves containing $S$
is not empty. 
We want to find a point $[C]\in M(k)$ corresponding to a curve  $C$.
For technical reasons  it is easier to concentrate on the case
when $[C]$ is a smooth point of $M$. Suppose that there is
 a smooth $k$-point $[C]$ and 
let $M_C\subset M$ denote  the unique irreducible component
containing $[C]$. Then $M_C$ is geometrically irreducible since
it has a smooth $k$-point.

The Galois group $\gal (\bar k/k)$ acts
on the irreducible components of $M_{\bar k}$ and it fixes $M_C$.
 Thus a first step to find $k$-rational points on $M$
is to find an irreducible subvariety 
$W\subset M_{\bar k}$ 
which is Galois invariant; that is,  defined over $k$ if $k$ is perfect.
(\ref{main.geom.thm}) asserts that this is always possible.

For an arbitrary field, for instance for $k=\q$,
finding $W$ helps very little in finding rational curves defined
over $k$. The situation is, however, much better over
finite fields, as we see in (\ref{geom=>ff.say}).
\end{say}

The main technical result of the paper
is a statement about SRC varieties over arbitrary perfect fields:

\begin{thm}\label{main.geom.thm}
Let $X\subset \p^N$ be a smooth, projective,
separably rationally connected variety 
of dimension $\geq 3$
over a  perfect field $k$. Let
 $S\subset X$ be a zero--dimensional  smooth subscheme.

Then  there is a nonempty locally closed subset $W$  of the Hilbert scheme
$\hilb(X)$ such that
\begin{enumerate}
\item $W$ parametrizes certain smooth rational curves $C\subset X$ 
of degree $d(X,S)$ which contain $S$.
\item $W$ is geometrically irreducible and smooth.
\item The restriction $T_X|_{C_w}(-S)$ is ample.

\item $d(X,S)$ is bounded from above 
 in terms of $\dim X$, $\deg X$ and $\deg S$. 
\end{enumerate}
\end {thm}

\begin{rem} In contrast with some similar results
(for instance \cite[3]{rcfg2}) we  no longer have
$k$-points in compactifications of $W$. We lose this
in Step 2 of (\ref{main.steps.say}).
\end{rem}

\begin{rem} Assume that $k$ is algebraically closed.
Let $R(X,S)\subset \hilb(X)$ denote the space of all rational
curves in $X$ passing through $S$. If $\alpha$ denotes a 
numerical equivalence class of curves, then $R(X)$ decomposes 
as a disjoint union $R(X,S)=\cup_{\alpha}R_{\alpha}(X,S)$.
It is reasonable to expect that for ``general'' $\alpha$,
the spaces $R_{\alpha}(X,S)$ are irreducible.
This would be a much stronger statement  than
(\ref{main.geom.thm}). There are, however, very few definitive
results in this direction.
The case when $X$ is homogeneous and $S=\emptyset$
is treated in \cite{kim-pandh}.
Lines on low degree hypersurfaces are discussed in
\cite{hmp}.
\end{rem}

The most straightforward application of (\ref{main.geom.thm})  is
to varieties over  pseudo-algebraically closed fields.
By definition, a
{\it  pseudo-algebraically closed field} or {\it PAC field}
is a field $k$ such that every geometrically irreducible
$k$-variety has a dense set of $k$-points. 
(See \cite{fried-jardin} for many examples and properties of
such fields.) 
(\ref{main.geom.thm}) now implies the following:

\begin{thm}\label{main.PAC.thm} 
Let $X$ be a smooth, projective,
separably rationally connected variety 
of dimension $\geq 3$
over a perfect PAC field $k$
and $S\subset X$ a zero--dimensional smooth subscheme.

Then there is a smooth, geometrically  rational curve 
$C_S\subset X$ which is defined over $k$ and contains $S$.

In particular, R-equivalence is trivial on $X(k)$
and  $\CH^0_0(X)=0$.\qed 
\end {thm}

The proof of the implication (\ref{main.PAC.thm}) $\Rightarrow$
(\ref{main.charp.thm})
 is  completely formal  using the fact
that the ultraproduct of infinitely many finite fields
whose orders go to infinity 
is PAC. See, for instance,  \cite{fried-jardin} for such techniques.

A more classical proof, given below, 
 relies on some basic boundedness results and
on the 
Lang-Weil theorem \cite{Lang-Weil}
(which is also at the core of the above ultraproduct result).

For a   quasi-projective scheme  $U\subset \p^N$ 
let  $\overline U$ denote  its closure
and $\partial U=\overline U\setminus U$  its boundary.
Define the  {\it degree} of a  reduced
projective scheme to be the sum of the degrees
of its irreducible components (even if they have different dimensions).
We call $N$, $\dim(U)$, $\deg(\overline U)$ and $\deg(\partial U)$
the {\it basic projective invariants} of a   reduced
 quasi-projective scheme 
$U$.

The set of all reduced
quasi-projective schemes with given  basic projective invariants
is bounded. That is, there is a scheme of finite type
$S$ and a locally closed subscheme $U\subset S\times \p^N$
such that every reduced
quasi-projective scheme with given  basic projective invariants
occurs among the fibers of $U\to S$.
(This is an easy consequence of the classical theorem, going back to Cayley,
that
the Chow variety of pure dimensional reduced closed subschemes of $\p^N$
of a given degree is bounded. See \cite{ho-pe} for a classical treatment
and \cite[I.3--4]{rcbook} for a modern update.)

Conversely, if $S$ is a scheme of finite type and
$U\subset S\times \p^N$ is a locally closed subscheme
then the basic projective invariants of the fibers form a bounded set.

Using basic existence results of the relative Hilbert scheme
(see, for instance, \cite[I.1]{rcbook}),
we get the following:

\begin{prop}\label{hilb.fin.lem}
 Let $X\subset \p^N$ be a projective variety
and $S\subset X$ a zero--dimensional subvariety.
Let $W_d\subset \hilb(X)$ be 
the locally closed  subset
of the Hilbert scheme parametrizing all
smooth curves of degree $d$ containing $S$.
Let $W_d^i\subset W_d$ be any irreducible component.

Then the basic projective invariants of $W_d^i$
can be bounded in terms of  $\dim X, \deg X,\deg S$ and $d$.\qed
\end{prop}

\begin{say}[Proof of (\ref{main.geom.thm}) $\Rightarrow$
(\ref{main.charp.thm})]\label{geom=>ff.say}

Let us now turn to the situation of (\ref{main.geom.thm}).
Let $\hilb(X,S)^s$ denote the Hilbert scheme of all smooth
rational curves  of $X$ which contain $S$.
By the deformation theory of Hilbert schemes containing a given subscheme
(see, for instance, \cite{mori-pn}), the assumption
 (\ref{main.geom.thm}.3) implies that  every point of $W$ 
is a smooth point of $\hilb(X,S)^s$.
(In our case $\hilb(X,S)^s$ can be identified with an open subset of the
Hilbert scheme of   $B_SX$, the blow up of $X$ at $S$, and the
more traditional study of $\hilb(B_SX)$ also gives this result,
cf.\  \cite[I.2.8]{rcbook}.)

Since $W$ is geometrically irreducible by assumption
(\ref{main.geom.thm}.2), this implies that 
 $W$ is contained in  a unique geometrically 
irreducible component $W'$ of  $\hilb(X,S)^s$.

By (\ref{hilb.fin.lem}), 
$W'$ is a quasi-projective variety in some $\p^N$, and
its basic projective invariants 
$N$, $\dim W'$, $\deg\overline{W'}$, $\deg\partial W'$ are bounded
 in terms of $\dim X$, $\deg X$  and $\deg S$.

The  Lang-Weil theorem  asserts that an $r$-dimensional
geometrically irreducible projective variety $W'$ has roughly $q^r$ points
defined over  $\f_q$, and the error term is
bounded by $C\cdot q^{r-1/2}$ 
with $C$  depending only on the basic projective
invariants of $W'$. The same remains true for quasi-projective
varieties, because the boundary has smaller dimension, hence its
points change only the error term.
(The boundary may not be geometrically irreducible,
but an {\em upper} bound for the number of points
holds without geometric irreducibility.)
 So one can easily find a lower
bound on $q$ in terms of the basic projective invariants, such that if
$q$ is larger than this bound then $W'$ has
$\f_q$-points.\qed
\end{say}

%%%%%%%%%%%%%%%%%%%%%%%%%%%%%%%%%%%%%%%%%%%%%%%%%%%%%%%%%%%%%%%%%%%%%%%%%%%
%%%%%%%%%%%%%%%%%%%%%%%%%%%%%%%%%%%%%%%%%%%%%%%%%%%%%%%%%%%%%%%%%%%%%%%%%%%

\begin{say}[Main steps of the proof of (\ref{main.geom.thm})]
\label{main.steps.say}

The construction of the geometrically irreducible
family of curves $W$ proceeds in four steps.

\begin{enumerate}
\item First, in lemma (\ref{through.one.point.lem}),
for every $P\in X(\bar k)$
 we construct geometrically irreducible
families (defined over $k(P)$) of maps passing through $P$
$$
V_P\subset \Hom (\p^1, X; (0{:}1)\mapsto P).
$$
These families were  introduced and used in \cite{rcloc}.
\item Then in (\ref{step.2.say}) 
 we find some auxiliary curves $C\subset X$
and  prove that the space of all  maps in $V_P$ 
whose image intersects $C$:
$$
V_{P,C}:=\{f\in V_P: f(1{:}0)\in  C\}
$$
is geometrically irreducible.
This relies on the Lefschetz--type results of \cite{rcfg1, rcfg2}
recalled in  (\ref{ratlef.lem}).
\item Third, starting with (\ref{step.3.say}),
 given $S=\{P_1,\dots,P_n\}$,
we construct a family $V_{S,C}$ 
of  reducible rational curves,
called combs, passing through $S$. Roughly speaking, over $\bar k$, 
$$
V_{S,C}:=V_{P_1,C}\times\cdots\times V_{P_n,C},
$$
and the curve corresponding to a point of $V_{S,C}$
consists of $C$ and  $n$ curves connecting $C$ to the points $P_1,\dots,P_n$.
(It is better to think of these as a curve mapping to $X$,
rather than a subscheme of $X$.)
This is the main new idea of the paper.
\item Finally, in (\ref{end.of.proof}),  we get  
the intended $W$ by smoothing the reducible curves in $V_{S,C}$.
\end{enumerate}
\end{say}

The easiest case, Step 1, is
handled in the following, essentially proved in  \cite[3.2]{rcloc}.
(See also 
\cite[44]{arko} or \cite[Thm.3,Rem.4]{rcfg2} where a much stronger
result is proved.)

\begin{lem}\label{through.one.point.lem}
Let $X\subset \p^m$ be a smooth, projective, SRC variety 
over a field  $k$ and $P\in X(k) $  a
point.
Then there is a family of rational curves
$$
F_P:V_P\times\p^1  \to  X,
$$ 
defined over $k$, 
with the following properties:
\begin{enumerate}
\item $F_P(V_P\times (0{:}1)) = P$.
\item $F_P:V_P\times (1{:}0)\to X$ is smooth, hence dominant.
\item $V_P$ is geometrically irreducible and smooth.
\item The pullback of $T_X$ to each curve $\p^1_v:=\{v\}\times \p^1$ is ample,
even after twisting by $\o(-3)$.
\item  The degree of the curves $F_P(\p^1_v)$ 
is bounded 
 in terms of $\dim X$ and $\deg X$. \qed
\end{enumerate}
\end{lem}

\begin{say}[Step 2] \label{step.2.say}
Next consider  the  space of rational curves
passing through a given point $P$ and meeting a given curve $C$. 
We start with the family of morphisms as in 
(\ref{through.one.point.lem}) and 
look at the evaluation map 
$$
f_P:
V_P  \to   X
$$
which sends $v\in V_P$ to $F_P(v, (1{:}0))$.
The inverse image
$f_P^{-1}(C)$  can be thought of as  a family  of 
 curves passing through $P$  and intersecting  $C$.
We would like $f_P^{-1}(C)$  to be geometrically irreducible.
We cannot provide a single curve $C$ that does the job.
For one thing, the curve should work simultaneously for many points
$P$, and for another, the curve must be defined over the base field.
On the other hand,  one can find a family of curves such that for each 
point  $P$ the general member of the family does work.
The precise   result we need is the following:

\begin{lem}\label{ratlef.lem} \cite[Thm.6,Cor.7,Rem.4]{rcfg2}
Let $X\subset \p^m$ be a smooth, projective variety over a field $k$.
Assume that $X$ is SRC. 
Then there is a nonempty family of rational curves, defined over $k$, 
$$
F:U\times \p^1\to X
$$
with the following properties:
\begin{enumerate}
\item
Let $Y$ be a geometrically irreducible variety and
$h:Y\to X$ a dominant morphism.
Then there is a dense open set $U^h\subset U$
such that for every  $u\in U^h(\bar k)$  
the fiber product
$Y\times_X \p^1_u$
is  irreducible.
\item $U$ is geometrically irreducible, smooth
and open in $\Hom(\p^1,X)$.
\item The pullback of $T_X$ to each curve $\p^1_u$ is ample.
\item  The degree of the curves $F(\p^1_u)$ 
is bounded  in terms of $\dim X$ and $\deg X$. 
\item  If $\dim X\geq 3$ then
$F$ is an embedding on every fiber $\p^1_u$.\qed
\end{enumerate}
\end{lem}
\end{say}

\begin{say}[Step 3]\label{step.3.say}
Now we move to  reducible curves passing through several points.
Assume that  our curves need to go
through the conjugation invariant set of  closed points $P_1,\dots, P_n$.
For each $P_i$ we choose a family $V_{P_i}$ 
defined over $k(P_i)$ as in 
(\ref{through.one.point.lem}) and also choose
another family as in (\ref{ratlef.lem}).
Out of these we get configurations of $n+1$ curves:
 an appropriate rational  curve $C$ and for each
$P_i$ a rational curve connecting $P_i$ to $C$.
\end{say}

\begin{defn} \label{defcomb}
A genus zero {\it comb} over $k$ with {\it $n$ teeth} is a 
reduced projective curve
of genus 0 (i.e. a curve $C$ with $h^1(C,\o_C)=0$)
having $n+1$ irreducible components over $\bar k$.
One component, defined over $k$, is called the {\it handle}.
The other $n$ components, $C_1,\dots, C_n$ 
are disjoint from each other and 
intersect  $C_0$
in $n$ distinct points. Every $C_i$ is smooth and rational.
The curves $C_1,\dots, C_n$ may not be 
individually defined
over $k$. 
A comb can be pictured as below:
$$
\begin{array}{c}
\begin{picture}(100,100)(40,-70)
\put(0,0){\line(1,0){180}}
\put(20,10){\line(0,-1){60}}
\put(40,10){\line(0,-1){60}}
\put(130,10){\line(0,-1){60}}
\put(160,10){\line(0,-1){60}}

\put(15,-60){$C_1$}
\put(35,-60){$C_2$}
\put(125,-60){$C_{n-1}$}
\put(155,-60){$C_n$}

\put(-25,-5){$C_0$}

\put(70,-30){$\cdots\cdots$}

\end{picture}\\
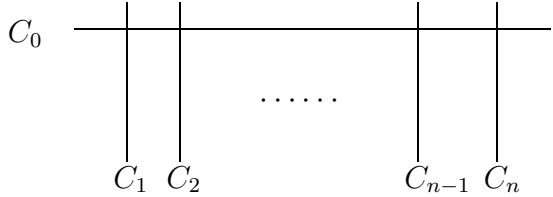
\mbox{Comb with $n$-teeth}
\end{array}
$$
\end{defn}

\begin{defn}\label{comb.defn}
Let $X$ be a smooth projective variety over  $k$ and 
$S\subset X$ a zero--dimensional smooth subset. Set $n=\deg S$.
A {\it $k$-comb through $S$} 
is
\begin{enumerate}
\item  a $k$-comb with $n$ teeth $C=C_0\cup\cdots\cup C_n$,
\item a $k$-morphism $f:C\to X$, and
\item a $k$-morphism $\sigma:S\to C\setminus C_0$ such 
that $f\circ \sigma=id_S$,
and each $C_i, i\geq 1$ contains exactly one point of $S_{\bar k}$.
\end{enumerate}

Set $S_{\bar k}=\{P_1,\dots,P_n\}$.
Over a perfect field, a comb through $S$ can also be given by
a $k$-morphism $g_0:\p^1\to X$ and a
conjugate set of $\bar k$-morphisms
$g_i:\p^1\to X$, $i=1,\dots, n$
such  that
\begin{enumerate}
\item  $g_i(0{:}1)=P_i$,  and 
\item  $g_i(1{:}0) = g_0(s(P_i))$ for $i=1,\dots, n$,
where $s:S\to \p^1$ is some embedding.
\end{enumerate}

A {\it family of combs through $S$} parametrized by a
scheme $T$ can be identified with   $n+1$ morphisms 
$$
G_i:T\times\p^1\to X \quad i=0,\dots, n
$$
with the expected properties.
\end{defn}

%%%%%%%%%%%%%%%%%%%%%%%%%%%%%%%%%%%%%%

\begin{say}[Construction of the combs]\label{comb.constr.say}

Let $X\subset \p^N$ be a smooth, projective,
separably rationally connected variety 
over a  perfect field $k$ and
 $S\subset X$  a zero--dimensional  smooth subscheme.
Set $S_{\bar k}=\{P_1,\dots,P_n\}$.

For a given $P_i\in S(\bar k)$ choose a family of maps
$$
F_{P_i}:V_{P_i}\times \p^1\to X,
$$
defined over $k(P_i)$,
as in  (\ref{through.one.point.lem}).
For  the conjugates $P_{i_1},\dots,P_{i_s}$  we choose  the conjugate families
$$
F_{P_{i_j}}:V_{P_{i_j}}\times \p^1\to X.
$$

Going through all the points of $S$,
 we have a conjugation invariant collection of $n$ families
$$
F_{P_i}:V_{P_i}\times \p^1\to X
$$
as in  (\ref{through.one.point.lem}).

Let $F:U\times \p^1\to X$ denote the family of curves obtained in
(\ref{ratlef.lem}).

For every $j$, let
$U(P_j)$ denote the subset of those
$$
(u,v_j)\in U\times V_{P_j}
$$
such that 
$$
F_{P_j}(v_j,(0{:}1))\in F(\p^1_u).
$$
The fiber of the projection 
$$
\pi_j: U(P_j)\to U
$$
over a point $u\in U$ consists of those
maps in the family $V_{P_j}$ which map $(0{:}1)$ to a point of 
$F(\p^1_u)$.
Thus by (\ref{ratlef.lem}) there is a nonempty open subset 
 $U^j\subset U_{\bar k}$
such that 
$\pi_j^{-1}(u)$ is geometrically irreducible for $u\in U^j(\bar k)$.
We can further assume that the $U^j$ are conjugates of each other.

As a $\bar k$-scheme, let
 $U(P_1,\dots, P_n)$ denote the subset of those
$$
(u,v_1,\dots, v_n)\in U\times \prod_{i=1}^nV_{P_i}
$$
such that 
$$
F_{P_i}(v_i,(0{:}1))\in F(\p^1_u) \quad  \forall\ i=1,\dots, n.
$$
By construction
$$
U(P_1,\dots, P_n)=U(P_1)\times_U\cdots\times_U U(P_n).
$$
Set   $U^S=U^1\cap \cdots\cap U^n\subset U$. Then the 
 fibers of the projection
$$
\Pi: U(P_1,\dots, P_n)\to U
$$
are geometrically irreducible  over $U^S$ since
$$
\Pi^{-1}(u)=\pi_1^{-1}(u)\times \cdots\times \pi_n^{-1}(u)
$$
and the product of geometrically irreducible schemes
is again geometrically irreducible.
Thus we conclude that:

\begin{claim} The variety $U(P_1,\dots, P_n)$ is defined over $k$ and it
has a unique geometrically irreducible component
$U^0(P_1,\dots, P_n)$ which dominates $U$.
\end{claim}

By our choices of the families $V_{P_i}$,
the construction of $U(P_1,\dots, P_n)$
is invariant under the action of $\gal(\bar k/k)$.
Thus we can henceforth  view $U(P_1,\dots, P_n)$ as a $k$-variety.
Not all of the points of $U^0(P_1,\dots, P_n)$ 
correspond to combs since the curves $C_i$ need not intersect
$C_0$ in distinct points.
In the construction, the intersection of $C_i$ and $C_0$
can be an arbitrary point of the set $C_0\cap F_{P_i}(V_{P_i}\times (1{:}0))$
which is dense in $C_0$ for general $C_0$.
 Therefore the combs
correspond to an open subset
$U^c(P_1,\dots, P_n)\subset U^0(P_1,\dots, P_n)$.
Combined with (\ref{through.one.point.lem}.4) and (\ref{ratlef.lem}.3)
this implies the following:

\begin{claim} $U^c(P_1,\dots, P_n)$
is a geometrically irreducible variety, defined over $k$.
It parametrizes combs through the subset $S$
$$
g:C=C_0\cup \cdots \cup C_n \to X
$$
with liftings $\sigma:S\to C$
 such that
$g^*T_X(-3\sigma(S))$ is ample on every irreducible component of $C$.
\end{claim}
\end{say}

\begin{say}[Step 4  of the proof of (\ref{main.geom.thm})]
\label{end.of.proof}

First we work in the space of genus zero stable curves 
$\bar M_0(X,S\into X)$ with base point set $S$ as
discussed in  \cite[sec.8]{arko}, and then we switch to the
Hilbert scheme of $X$.

By the universal property of $\bar M_0(X,S\into X)$, there is a
moduli map 
$$
m:U^c(P_1,\dots, P_n)\to \bar M_0(X,S\into X).
$$
We claim that there is  a  unique irreducible
component    $U_1\subset \bar M_0(X,S\into X)$ containing 
 $m(U^c(P_1,\dots, P_n))$.

Since  $g^*T_X(-3\sigma(S))$ is ample on every irreducible component of $C$,
an easy lemma (cf. \cite[18]{arko}) implies that 
$H^1(C, g^*T_X(-\sigma(S)))=0$.  Thus by \cite[42.1]{arko}
 the space $\bar M_0(X,S\into X)$ has a unique irreducible
component  containing the image of any point of
 $U^c(P_1,\dots, P_n)$. The latter is
irreducible, so this must be the same component, call it $U_1$,  for
every point.

Over $\bar k$, a general point of $U_1$
corresponds to an irreducible curve  
 $g_u:D_u\cong \p^1\to X$ with $\sigma_u:S\to D_u$ by \cite[42.3]{arko}.
By semicontinuity, $H^1(D_u, g_u^*T_X(-3\sigma_u(S)))=0$
for $u$ in an open subset of $U_1$. The pull back
$g_u^*T_X(-\sigma_u(S))$ is a vector bundle on $\p^1$,
hence a direct sum of line bundles.
The above vanishing 
 implies that every direct summand of $g_u^*T_X(-\sigma_u(S))$
has degree at least $2\deg S-1\geq 1$.
If $\dim X\geq 3$ then  by \cite[II.3.14.3]{rcbook} this implies that
there is an open subset $U_2\subset U_1$ parametrizing
embeddings $\p^1\into X$ passing through $S$.
(If $\dim X=2$, we get  immersions by \cite[II.3.14.2]{rcbook}.)

Pick any $u\in U_2(\bar k)$. Since $g_u:D_u\into X$ is an
embedding, there are no automorphisms of $D_u$ commuting with
$g_u$. Thus over $U_2$ we have a fine moduli problem and 
there is  a universal family ${\mathbf D}_2\to  U_2$,
 see \cite[Thm.2.iii]{FuPa}.
The graph of the universal morphism 
$$
{\mathbf g}:{\mathbf D}_2\into U_2\times X
$$
defines an embedding $U_2\to \hilb(X)$ with image $W$.
\qed

\end{say}

\begin{say}[Proof of (\ref{finite-R.thm})]

Let $X$ be a smooth, projective,
separably rationally connected variety over the finite field $k$. 
Pick  $P,Q\in X(k)$. Set $S=P\cup Q$.
By (\ref{main.charp.thm}), if $\dim X\geq 3$ and 
if $|k|>\Phi(\deg X ,\dim X ,\deg S)$, 
then there is a smooth,  rational curve
$C_S\subset X$ containing $S$,
hence $P$ and $Q$ are R-equivalent.
If $\dim X=2$ then we work with $X\times \p^1$.

Set
$\Psi(a,b)=\Phi(a,b,2)$ for $b\geq 3$ and $\Psi(a,2)=\Phi(3a,2,2)$.
($3\deg X$ is the degree of $X\times \p^1$ under the Segre embedding
when $X$ is a  surface.)
\qed
\end{say}

\begin{say}[Proof of (\ref{local-R.thm})]

Let $X$ be a smooth, projective,
separably rationally connected variety over the local field $K$. 
Pick $P,Q\in X(K)$. Set $S=P\cup Q$.
The base field $K$ is the quotient field of a 
complete valuation ring $R$ whose
residue field is $k$. By assumption there is a model 
$X_R\to\spec R$ of $X$ 
such that the special fibre $X_k$
satisfies the conditions of (\ref{finite-R.thm}).
Since $X_R$ is proper and $R$ is a valuation ring, the points 
$P,Q$ give rise to sections $P_R,Q_R$ of $X_R\to\spec R$,
hence they have well-defined reductions $P_k,Q_k\in X_k$.

Assume first that $P_k\neq Q_k$.
By (\ref{main.charp.thm}), 
if $|k|>\Phi(\deg X ,\dim X ,\deg S)$, 
then there is a smooth, geometrically  rational curve
$C_S\subset X_k$ containing $S_k$.
Moreover, the twisted restriction $T_X|_{C_{S_k}}(-S_k)$ is ample.

We apply \cite[51]{arko} to $X_R\to\spec R$
with base points $P_R\cup Q_R\into X_R$.
(The current  $R$ corresponds to  $S$ of \cite{arko},
the $P\to S$ there is 2 copies of $S$ and $t:P\to X_R$
is the imbedding whose image is $P_R\cup Q_R$.)
By  \cite[51.5]{arko}
we obtain a 
smooth $R$-scheme $U\to R$ 
with a point $u\in U(k)$ and a family of rational curves
$U\times \p^1\to X_R$ passing through $P_R$ and $Q_R$.
By the Hensel lemma, $\phi_k:\spec k\to u$ lifts to
$\phi_R:R\to U$.
The image of  $\phi_R(\spec K)$ corresponds to a 
 rational curve $C\subset  X_K$ going through $P$ and $Q$.
Thus $P$ and $Q$ are R-equivalent.

If $P_k=Q_k$ then we work with $X\times \p^1$
and consider the points $P\times (0{:}1)$ and $Q\times (1{:}0)$.
These have reductions  $P_k\times (0{:}1)$ and $Q_k\times (1{:}0)$
which are now different. The previous argument
 gives a rational curve on $X\times \p^1$
passing through $P\times (0{:}1)$ and $Q\times (1{:}0)$.
This again implies that $P$ and $Q$ are R-equivalent.

\qed

\end{say}

%%%%%%%%%%%%%%%%%%%%%%%%%%%%%%%%%%%%%%%%%%%%%%%%%%%%%%%%%%%%%%%%%%%%%%%%%%

\begin{say}[Proof of (\ref{finite-CH.thm})]

Let $X$ be a smooth, projective,
separably rationally connected variety over the finite field $k$. 
Let $Z, Z'$ be two effective zero cycles
of the same degree  on $X$. Let $S$ be the support of  $Z+Z'$.

By (\ref{main.charp.thm}), 
if $|L|>\Phi(\deg X ,\dim X ,\deg S)$, 
then there is a smooth  rational curve
$C_L\subset X_L$ containing $S$.
Thus $Z-Z'$ is trivial in $\CH_0(X_L)$.

The pull back  followed by the push forward
$$
\CH_0(X_k)\to \CH_0(X_L)\to \CH_0(X_k)
$$
is multiplication by  $\deg (L/k)$.
This proves that $Z-Z'\in \CH_0(X)$ is torsion whose order
divides $\deg (L/k)$.

A finite field has extensions of arbitrary degree.
Taking two such extensions of  degrees
$a$ and $a+1$ for $a>\Phi(\deg X ,\dim X ,\deg S) $,
we conclude that 
$Z-Z'\in \CH_0(X)$ is torsion whose order
divides $\gcd(a,a+1)=1$.\qed
\end{say}

\begin{say}[Proof of (\ref{local-CH.thm})]\label{ct.proof}

Let $X$ be a smooth, projective,
separably rationally connected variety over the local field $K$. 

Let $X_R\to R$ be a smooth projective model of $X_K$
with special fiber $X_k$. Every finite extension
$k'\supset k$ corresponds to a unique unramified
extension $K'\supset K$ 
(This is essentially
the theory of Witt rings as explained in
\cite[Chap.II]{serre-loc}.)

Thus there are two extensions $K_i/K$ of relatively  prime
degree such that $X(K_i)\neq \emptyset$. 
As in the previous proof, it is enough to prove
that $\CH_0^0(X_{K_i})=0$.
Changing notation, it is enough to prove
that $\CH_0^0(X)=0$ if there is a point $P\in X(K)$.

Let $Q\in X$ be any closed point of degree $d$.
It is enough to prove that $Q-dP\in  \CH_0^0(X)$ is zero.
Set $L=K(Q)$ and let $Q_1\in X(L)$ be a point over  $Q$.
We have proved that $Q_1$ and $P$ are R-equivalent
in $X_L$, thus $Q_1-P$ is zero in  $\CH_0^0(X_L)$.
Hence its push forward $Q-dP$ is zero in  $\CH_0^0(X)$.
\qed  

\end{say}

\begin{ack}   We   thank J.-L.\ Colliot-Th\'el\`ene for
many useful conversations, questions and a long list of corrections
and improvements.
The above simple  proof in (\ref{ct.proof}) is due to him.
Partial financial support
to the first author 
 was provided by  the NSF under grant numbers 
DMS-9970855 and DMS02-00883. Partial financial support
to the second author 
 was provided by  OTKA under grant numbers T029525 and T31984.
The paper was written while the first author
visited the  Mathematical Institute of the Hungarian Academy of Sciences.
 \end{ack}

\vskip1cm

\noindent Princeton University, Princeton NJ 08544-1000

\begin{verbatim}kollar@math.princeton.edu\end{verbatim}

\vskip0.5cm

\noindent Mathematical Institute of the Hungarian Academy of Sciences,

\noindent 1364 Budapest, PO Box 127

\begin{verbatim}endre@renyi.hu\end{verbatim}

\end{document}